\documentclass[a4paper,11pt]{amsart}

\usepackage[english]{babel}
\usepackage{amsmath, enumerate, amsfonts, amssymb, amsthm}

\newtheorem{defin}{\bf Def\mbox{}inition}[subsection]
\newtheorem{theo}[defin]{\bf Theorem}
\newtheorem{prop}[defin]{\bf Proposition}
\newtheorem{lem}[defin]{\bf Lemma}
\newtheorem{cor}[defin]{\bf Corollary}

\newtheorem{rem}[defin]{\bf Remark}

\newtheorem*{clai*}{\bf Claim}
\newtheorem*{rem*}{\bf Remark}
\newtheorem*{nota*}{\bf Notation}
\newtheorem*{prop*}{\bf Proposition}

\newtheorem{exI}{\bf Example}

\newcommand{\dps}{\displaystyle}

\newcommand{\R}{\mathbb{R}}
\newcommand{\C}{\mathbb{C}}
\newcommand{\N}{\mathbb{N}}

\renewcommand{\k}{\mathbf{k}}

\newcommand{\ring}{\mathcal{R}}

\newcommand{\CC}{\mathcal{C}}
\newcommand{\FF}{\mathcal{F}}
\newcommand{\Frac}{\mathrm{Frac}}

\newcommand{\QQ}{\mathcal{Q}}
\newcommand{\modQ}{{\mathrm{mod}\QQ}}
\newcommand{\PP}{\mathcal{P}}

\newcommand{\GG}{\mathcal{G}}

\newcommand{\spec}{\mathrm{Spec}}
\newcommand{\specm}{\mathrm{Specm}}

\renewcommand{\O}{\mathcal{O}} 

\newcommand{\FD}{\hat{\mathcal{D}}} 
\newcommand{\FDn}{\hat{\mathcal{D}}_n}

\newcommand{\dx}[1]{\partial _{x_{#1}}}

\newcommand{\ddx}{\partial _x}

\newcommand{\ord}{\mathrm{ord}} 
\newcommand{\gr}{\mathrm{gr}} 
\newcommand{\lm}{\mathrm{lm}} 
\newcommand{\lc}{\mathrm{lc}} 
\newcommand{\lt}{\mathrm{lt}} 
\newcommand{\ND}{\mathcal{N}} 
\newcommand{\Exp}{\mathrm{Exp}} 

\title{Generic and comprehensive standard bases}

\author{Rouchdi BAHLOUL} \address{Department of Mathematics, Faculty
of Science, Kobe University, 1-1, Rokkodai, Nada-ku, Kobe 657-8501,
Japan} \email{rouchdi@math.kobe-u.ac.jp}

\begin{document}

\begin{abstract}
Parametric Gr\"obner bases have been studied for more than 15 years and
are now a further developed subject. Here we propose a general study of
parametric standard bases, that is with local orders. We mainly focus on
the commutative case but we also treat the case of differential operators
rings. We will be concerned by two aspects: a theoretical aspect with
existence theorems and a practical aspect devoted to how we can explicitely
compute such objects when the given data are algebraic.
We believe that parametric standard bases are important for both aspects.
From a theoretical point of view, they constitute a strong tool for proving
constructive results. From a practical one, they provide a tool for studying
explicitely local objects associated with parametric algebraic ideals.
\end{abstract}

\maketitle


\section*{Introduction}\label{sec:intro}

Parametric Gr\"obner bases have more than 15 years old and have been
studied or used by several authors: P.~Gianni \cite{gianni},
M.~Lejeune-Jalabert and A.~Philippe \cite{lej-phi}, D.~Bayer et al.
\cite{bgs}, A.~Assi \cite{assi}, T.~Becker \cite{becker}, M.~Kalkbrener
\cite{kalkbrener}, E.~Fortuna et al. \cite{fgt}, A.~Montes \cite{montes},
V.~Weispfenning \cite{weisp92} for the notion of comprehensive Gr\"obner
basis and more recently \cite{weisp03} with a canonical treatment; see also
an alternative approach by T.~Sato and A.~Suzuki \cite{satosuzuki}, etc.\\ 
All these constructions take place in polynomial rings with well-orderings.
In 1997, T.~Oaku \cite{oaku} used parametric Gr\"obner bases in rings of
algebraic differential operators. This work inspired A.~Leykin \cite{leykin}
and U.~Walther \cite{walther}. Again, these constructions concerned
well-orderings.

Now, concerning the local situation, in a recent paper A.~Fr\"uhbis-Kr\"uger
\cite{fruhbis} made use of a parametric approach to standard bases for the
study of families of singularities. However, it seems that there does not
exist any general study concerning standard bases with parameters.

The motivation of this paper is a natural question. We know how to
construct parametric Gr\"obner bases for ideals in\footnote{Throughout
the paper, the symbol $\k$ denotes a field} $\k[a,x]=
\k[a_1,\ldots, a_m, x_1, \ldots, x_n]$ (we see $x$ as the main variables
and $a$ as a system of parameters and $\k$ is a field) and how to study
their behaviour when we specialize the parameter $a=c$ with $c\in \k^m$.
What can we say about an ideal in $\C\{a,x\}$ when we specialize the
variable $a=c$ with $c\in \C^m$ in a neighbourhood of $0$? Moreover
for an ideal in $\k[a,x]$ and a local (or arbitrary) order on the monomials
in $x$, can we make explicite calculations?

These are natural questions, the first one being theoretical (in the
sense that we cannot have finite algorithms) and the second one being
practical.

The second question arises for example if we want to study the behaviour
of the germ at $0$ of the variety $V(I_{|a=c}) \subset \k^n$ when $c$
runs over $\k^m$ for $I\subset \k[x,a]$.

The purpose of this article is to answer these two qustions. We will also
treat the case of rings of differential operators.
However the proofs being the same, we felt it would be more convenient to
the reader to give the statements and the proofs in the commutative case
and separately to give only the statements in the non commutative case.

In order to motivate the reading, let us examine a trivial but instructive
example.

\begin{exI}\label{exemple}
Let $f=f(a,x_1,x_2)=a x_2-x_1 x_2+ x_1 \in \k[a][x_1,x_2]$. Let $\prec$
be a local order on the terms $x_1^i x_2^j$ such that the leading term of
$f$ is $x_2$ (here $a$ is one parameter). Let $I= \k[a][x_1, x_2] \cdot f$
and $\hat{I}=\k[a][[x_1, x_2]] \cdot f$. In this trivial situation,
$f(c,x_1, x_2)$ is a standard basis of $I_{|a=c}$ and $\hat{I}_{|a=c}$ for
any $c\in \k$ where the leading term is $x_2$ for any $c\ne 0$ and is
$x_1$ for $c=0$. However, this standard basis is not reduced. Let us
(formally) reduce $f$ as if $a$ were a constant. Rigorously, we work in
$\Frac(\k[a])[[x_1, x_2]]$ where $\Frac$ denotes the fraction field.
The reduction of $f$ is
\[x_2 + x_1/a + x_1^2/a^2 + x_1^3/a^3 + \cdots .\]

As we can see, this is neither in the ring $\k[a][x_1,x_2]$ nor in
$\k[a][[x_1,x_2]]$ but in $\k[a][a^{-1}][[x_1,x_2]]$, that is power series
with coefficients in the localization of $\k[a]$ with respect to $a$,
which can be seen as a subring of $\Frac(\k[a])[[x_1,x_2]]$.
\end{exI}

In order to justify the construction adopted in this paper, let us
make a remark: as we saw here, if we don't ask our ``generic standard
basis'' (this term shall be defined later) to be reduced then it can
be chosen in the ring where the given generators are (this fact
is general as shall be proved in section 2) but if we want to
construct a ``generic reduced standard basis'' then we will have to
work in some localization $\CC[h^{-1}][[x]]$ with $h\in \CC$
where $\CC$ is the ring of parameters, here above we had
$\CC=\k[a]$. Therefore, the natural environment in which our objects
are is the ring $\Frac(\CC)[[x]]$.

Moreover, let us say that this necessity of localizing which (as been
said) occures when we need reduced generic standard bases, shall be
the main difference with the global case (that is the case of
polynomial rings with global orders).

Let us give the contents of the paper. In the first section, we give
some recalls without the proofs. This concerns division theorems and
standard bases in $\k[[x]]$ and $\C\{x\}$. We will introduce the notion
of truncated division (it appeared after discussing with A.~Assi that
this notion has been already introduced in \cite{assiT}).

Section 2 constitutes the heart of the paper with the notion of
generic standard basis and generic reduced standard basis (for short
gen.s.b and gen.red.s.b).

In section 3, we show how we can in an algorithmic way construct a
gen.s.b if we start with an ideal in $\k[a,x]$. Let us point out
the fact that in this section, we shall work with an arbitrary order
(not necessarily a local order).

Section 4 contains direct applications.
The first one is the existence of a comprehensive standard basis
for an ideal in $\C\{a,x\}$ by using that of a gen.s.b The second
application gives an example of a possible use of gen.s.b: it concerns
the Hilbert-Samuel polynomial of the germ of variety associated with
an ideal in $\C\{a,x\}$.

In Section 5, we extend our result to ideals in rings of differential
operators. We also treat the homogenized version.
As said, proofs are not given since they are the same as in the
commutative case.\\

This paper has been announced in \cite{B-caalias} (here some changes have
been made with more unified definitions).
A preliminary work was the preprint \cite{B-genGF} where we applied
gen.red.s.b to study the local Gr\"obner fan of an analytic ideal in rings
of differential operators. We plan to apply this work to a ``parametric''
study of differential systems.

In order to keep a reasonable size for the present paper, we restricted
ourselves to direct applications or illustrations. However a more substantial
use of generic standard bases (in rings of differential operators) is
made in \cite{B-polygen} where we study the local Bernstein-Sato polynomial
associated with a deformation of a singularity.\\
{\bf Acknowledgements.}
I would like to thank A. Assi, M. Granger and N. Takayama. Thanks to
their advise, remarks or questions, my work on parametric standard
bases has been considerably improved. Thanks go to O.~M.~Abderrahmane
for pointing out to me the reference \cite{broughton}.\\
This work is made under the support of the FY2003 JSPS Postdoctoral
Fellowship.

\section{Recalls on divisions and standard bases}

In the following subsections, we will give some recalls. For the
proofs, the reader can refer to \cite{cg}. We also introduce truncated
divisions as in (Assi, \cite{assiT}),  they will play an important role
in the sequel.

\subsection{Division theorem}\label{subsec:div}

Let $\prec$ be an order on the terms $x^\alpha=x_1^{\alpha_1}
\cdots x_1^{\alpha_n}$, $\alpha\in \N^n$. If it is compatible with products,
we say it is a \emph{monomial order}. An order $\prec$ is a \emph{well
order} if any set of monomials have a minimum. A monomial order
is a well order if and only if: $x^\alpha \succeq 1$, for any $\alpha \in
\N^n$ (we also say `global order'). On the opposite side, a \emph{local
order} is a monomial order such that $x^\alpha \preceq 1$, $\forall \alpha
\in \N^n$. Note that we will denote by the same symbol $\prec$ the
order induced on $\N^n$.

Let $\ring$ be one of the following: $\k[x]$, $\C\{x\}$, $\k[[x]]$.
Let $\prec$ be a global order in the first case, otherwise it is a local
order.
Let $f$ be non zero in $\ring$. It has a unique writing as $f=\sum_{\alpha
\in \N^n} c_\alpha x^\alpha$ where $c_\alpha$ is in $\k$ or $\C$ and
this sum is finite if $f\in \k[x]$.

Define the Newton diagram $\ND(f)$ of $f$ as the set of $\alpha \in
\N^n$ such that $c_\alpha \ne 0$. Define the leading exponent of $f$
(w.r.t. $\prec$) as $\exp_\prec(f)=\mathrm{max}_\prec \ND(f)$. If there
is no confusion, we omit the subscript $\prec$. We then define the
leading term, leading coefficient and leading monomial of $f$ as:
$\lt(f)=x^{\exp(f)}$, $\lc(f)=c_{\exp(f)}$, $\lm(f)= \lc(f) \lt(f)$.

Let us now recall the division theorem with unique quotients and
remainder as in \cite{cg}.
For $e_1,\dots, e_r$ in $\N^n$, consider the following partition of $\N^n$.
\begin{itemize}
\item
$\Delta_1=e_1+\N^n$ and for $j=2,\ldots, r$, $\Delta_j=(e_j+\N^n)
\smallsetminus (\Delta_1 \cup \cdots \cup \Delta_{j-1})$.
\item
$\bar{\Delta}=\N^n \smallsetminus (\Delta_1 \cup \cdots \cup
\Delta_r)$.
\end{itemize}

\begin{theo}[\cite{cg}]\label{theo:div}
Let $g_1,\ldots,g_r$ be non zero elements in $\ring$ and consider the
partition associated with the $\exp(g_j)$. For any $f\in \ring$, there
exists a unique $(q_1,\ldots, q_r,R) \in \ring^{r+1}$ such that:
\begin{itemize}
\item[(i)] $f=q_1 g_1 + \cdots + q_r g_r+R$
\item[(ii)] for any $j$, $q_j=0$ or $\ND(q_j)+\exp(g_j) \subset
\Delta_j$
\item[(iii)] $R=0$ or $\ND(R) \subset \bar{\Delta}$.
\end{itemize}
$R$ is called the remainder of the division of $f$ by $g_1,\ldots,
g_r$ (w.r.t. $\prec$).
\end{theo}
As a consequence, we have:
\begin{equation}\label{eq:maxexp}
\exp(f)=\mathrm{max}_\prec \{\exp(q_1 g_1), \ldots, \exp(q_r g_r),
\exp(R)\} \textrm{ and } \exp(R) \in \bar{\Delta}.
\end{equation}

For the proof in the global case, see (\cite{cg}, Th. 1.5.1). Suppose $\prec$
is local and let us show how we can recover this theorem from \cite{cg}.
In (\cite{cg}, Th. 1.5.1) the authors proved the same result but for an
order $<_w$ defined in this way: they fix  a weight vector $w\in \R^n$ with
$w_i \le 0$ and a well order $<_0$ and they define $<_w$ in a lexicographical
way by $w$ and the inverse of $<_0$. By Robbiano's theorem \cite{robbiano},
it is easy to see that our order $\prec$ is of this form. Thus we are under
the hypothesis of (\cite{cg}, Th 1.5.1).

Let us recall the main steps of the proof of the theorem, this will be
useful in the sequel.\\
$\bullet$ Put $(f^{(0)}, q_1^{(0)}, \ldots, q_r^{(0)}, R^{(0)})=
(f, 0, \ldots, 0)$.\\
$\bullet$ For $i\ge 0$, if $f^{(i)}=0$ then put

$(f^{(i+1)}, q_1^{(i+1)}, \ldots, q_r^{(i+1)}, R^{(i+1)}) =
(f^{(i)}, q_1^{(i)}, \ldots, q_r^{(i)}, R^{(i)})$.\\
$\bullet$ If $\exp(f^{(i)}) \in \bar{\Delta}$ then

$(f^{(i+1)}, q_1^{(i+1)}, \ldots, q_r^{(i+1)}, R^{(i+1)}) = (f^{(i)}-
\lm(f^{(i)}), q_1^{(i)}, \ldots, q_r^{(i)}, R^{(i)} +\lm(f^{(i)}))$.\\
$\bullet$ If not, then let $j=\min\{k\in \{1,\ldots, r\} , \exp(f^{(i)})
\in \Delta_k\}$ and put

$\dps f^{(i+1)}= f^{(i)} - \frac{\lc(f^{(i)})}{\lc(g_j)} \cdot
x^{\exp(f^{(i)})-\exp(g_j)} \cdot g_j$,

$q_j^{(i+1)}= q_j^{(i)}+ \frac{\lc(f^{(i)})}{\lc(g_j)} \cdot
x^{\exp(f^{(i)})-\exp(g_j)}$,

for $l\ne j$, $q_l^{(i+1)}=q_l^{(i)}$ and $R^{(i+1)}=R^{(i)}$.

In this construction, we have $(r+2)$ sequences $f^{(i)}$, $q_1^{(i)}
\ldots, q_r^{(i)}$ and $R^{(i)}$ with $i \in \N$ satisfying $f=\sum_j
q_j^{(i)} g_j +R^{(i)} +f^{(i)}$. The first step consists in showing
that these sequences converge in $\k[[x]]$ for the
$(x_1,\ldots,x_n)$-adic topology (in particular, the limit of
$f^{(i)}$ is zero). The second step which is much harder is to prove
that if the inputs are in $\C\{x\}$ then so are the outputs.

As an easy consequence, we have:
\begin{lem}\label{lem:denom}
Let $\CC$ be a commutative integral ring and $\FF$ a field containing $\CC$.
Let $f, g_1,\ldots,g_r$ be in $\CC[[x]]$. Let us
consider the division of $f$ by the $g_j$'s in $\FF[[x]]$
w.r.t. $\prec$: $f=\sum q_j g_j +R$. Then the coefficients of $R$ and
of the $q_j$'s have the following form:
\[\frac{c}{\prod_{j=1}^r \lc(g_j)^{d_j}} \text{ where } c \in \CC, \,
d_j \in \N.\]
\end{lem}

Now let us end this subsection with the notion of \emph{truncated
division}.

\begin{defin}\label{def:truncated}
Given $f,g_1,\ldots,g_r \in \ring$. Let $f=\sum_j q_j g_j +R$ be the division
as in the theorem above. We define the truncated division of $f$ by
$g_1,\ldots,g_r$ to be:\\
$\bullet$ if $R=0$: $f=\sum_j q_j g_j$.\\
$\bullet$ if $R\ne 0$: $f=\sum_j q_j^{(i_0)} g_j + f^{(i_0)}$ where $i_0$
is the minimal $i$ such that $\exp(f^{(i)}) \in \bar{\Delta}$.

We call $f^{(i_0)}$ the remainder of the truncated division of $f$
by the $g_j$.
\end{defin}

\begin{rem}\label{rem:algebraic}
\noindent
\begin{itemize}
\item[(1)]
For this division, properties (ii) and (iii) of the division thereom are not
satisfied in general, however, the relation (\ref{eq:maxexp}) is satisfied.
\item[(2)]
In the second case of the definition, we have $q_j^{(i_0)} \in k[x]$ thus
\[f^{(i_0)} \in f + \sum_j \k[x] \cdot g_j.\]
\end{itemize}
\end{rem}

\subsection{Standard bases}

Here again, for the proofs, we refer to \cite{cg}.

We still denote by $\ring$ one of the rings $\C\{x\}$, $\k[[x]]$. Let
$J$ be an ideal in $\ring$ and consider the set of leading exponents
of $J$ (w.r.t. $\prec$):
\[\Exp(J)=\big\{ \exp(f); f \in J \smallsetminus \{0\} \big\}.\] 
This is a subset of $\N^n$ stable by sums, thus by the usual Dickson
lemma, we have:

\begin{defin}
There exists $G=\{g_1,\ldots, g_r\} \subset J$ such that
\[\Exp(J)=\bigcup_{j=1}^r (\exp(g_j) + \N^n).\]
Such a set $G$ is called a standard basis of $J$ (w.r.t. $\prec$).
\end{defin}

As a consequence of the division theorem, the following holds.
\begin{prop}[\cite{cg}, Cor. 1.5.4]\label{prop:rest0}
Let $G=\{g_1,\ldots,g_r\}$ in $J$, these two statements are equivalent:\\
$\bullet$ For any $f\in J$,  the remainder of the division of $f$ by
$G$ is zero.\\
$\bullet$ The set $G$ is a standard basis of $J$.
\end{prop}

Now we define the {\bf $S$-function} of $f$ and $g$ in $\ring$ as
$S(f,g)=\lc(g) m f- \lc(f) m' g$ where $m=m_0/\lt(f)$, $m'=m_0/\lt(g)$
and $m_0=\mathrm{lcm}(\lt(f), \lt(g))$.

As in the polynomial case with well orders \cite{buchberger}, we have a
Buchberger type criterion:
\begin{prop}[\cite{cg}, Prop. 1.6.2]\label{prop:Scriterion}
Let $G=\{g_1,\ldots,g_r\}$ be a set of generators of an ideal $J
\subset \ring$. Then $G$ is a standard basis of $J$ if and only if:
for any $j,j'$, the remainder of the division of $S(g_j, g_{j'})$ by
$G$ is zero.
\end{prop}

As a consequence, one can construct a standard basis by using the
Buchberger algorithm \cite{buchberger}, which consists, starting from a
set of generators $G_0$, in adding all the non zero remainders of the
division of $S(g,g')$ by $G_0$, with $g,g' \in G_0$ and, calling $G_1$
this new set, continuing the same process with $G_1$, etc. If we
denote $E(G)=\cup_g(\exp(g) +\N^n)$ then the termination of the
algorithm is assured by the fact that $E(G_0) \subset E(G_1) \subset
\cdots$ and by noetherianity of $\N^n$ and finally by
Prop.~\ref{prop:Scriterion}.

\begin{rem}[truncations]
In Buchberger algorithm, we can use truncated divisions as well.
The process will stop for the same reasons.
\end{rem}

By using truncated divisions in Buchberger algorithm we have the following
easy consequence (see Remark \ref{rem:algebraic}):

\begin{lem}\label{lem:algebraicSB}
Let $\ring_0$ be a subring of $\k[[x]]$ such that $\k[x] \subset \ring_0$.
Let $f_1, \ldots, f_q \in \ring_0$. Let $I$ (resp. $\hat{I}$) be the ideal
of $\ring_0$ (resp. of $\k[[x]]$) generated by the $f_j$ then there exists a
standard basis $G$ of $\hat{I}$ in
\[\k[x] \cdot f_1 + \cdots + \k[x]\cdot f_q.\]
In particular, $G$ is in $I$.
\end{lem}

Let us end these preliminaries with the notion of reduced standard basis.

\begin{defin}
A standard basis $G=\{g_1,\ldots,g_r\}$ of $J\subset \ring$ is said to
be
\begin{itemize}
\item
minimal if for any $F\subset \N^n$, we have:

$\Exp_{<_L^h}(J)=\bigcup_{e\in F} (e+\N^{2n+1}) \Rightarrow
\{\exp(g_1),\ldots,\exp(g_r)\} \subseteq F$.
\item
reduced if it is minimal and if for any $j=1,\ldots,r$, $\lc(g_j)=1$
and

$(\ND(g_j)\smallsetminus \exp(g_j)) \subset (\N^n \smallsetminus
\Exp(J))$.
\end{itemize}
\end{defin}

\begin{lem}\label{lem:reducedSB}
Given an ideal $J \subset \ring$ and a local
order $\prec$, a reduced standard basis exists and is unique.
\end{lem}
\begin{proof}
The unicity is left to the reader. Let us sketch the existence. Let
$G_0$ be any standard basis. By removing unecessary elements we may
assume $G_0$ to be minimal. Set $G_0=\{g_j; 1\le j \le r\}$. For
any $j$, divide $g_j-\lm(g_j)$ by $G_0$ and denote by $r_j$ the
remainder. The set $\{(\lm(g_j)+r_j)/\lc(g_j); \, 1\le j\le r\}$ is
then the reduced standard basis of $J$.
\end{proof}

\section{Generic standard bases}

Let $\CC$ be a commutative integral unitary ring for which we denote
by $\FF$ the fraction field, by $\spec(\CC)$ the spectrum and by
$\specm(\CC)$ the maximal spectrum.
For any ideal $\mathcal{I}$ in $\CC$, we denote by
$V(\mathcal{I})=\{\PP \in \spec(\CC); \mathcal{I} \subset \PP\}$ the
zero set defined by $\mathcal{I}$.

For any $\PP$ in $\spec(\CC)$ and $c$ in $\CC$, denote by $[c]_\PP$
the class of $c$ in $\CC/\PP$ and by $(c)_\PP$ this class viewed in the
fraction field $\FF(\PP)=\Frac(\CC/\PP)$. The element $(c)_\PP$ is
called the {\bf specialization} of $c$ {\bf into} $\PP$.

We naturally extend these notations to elements in $\CC[[x]]$ and we
extend $(\cdot )_\PP$ to elements of $\FF[[x]]$ for which the
denominators of the coefficients are in $\CC \smallsetminus \PP$, i.e.
$\CC_\PP[[x]]$ where $\CC_\PP$ is the localization w.r.t. $\PP$.

Now, given an ideal $J \subset \CC[[x]]$, we define the specialization
$(J)_\PP$ of $J$ into $\PP$ as the ideal of $\FF(\PP)[[x]]$ generated
by all the $(f)_\PP$ with $f\in J$.

\subsection{Generic standard basis on an irreducible affine scheme}

\

Fix a prime ideal $\QQ$ in $\CC$.
Let us start with some notations. We denote by $\QQ[[x]]$ the ideal of
$\CC[[x]]$ made of elements with all their coefficients in $\QQ$.
For $h\in \CC$, we denote by $\CC[h^{-1}]$ the localization of $\CC$ w.r.t.
$h$. The ring $\CC[h^{-1}][[x]]$ shall be seen as the subring of $\FF[[x]]$
made of elements with coefficients $\frac{c}{c'}$ such that $c'$ is a power
of $h$. In the latter, if all the $c$ are in $\QQ$, we obtain an ideal
denoted by $\QQ[h^{-1}][[x]]$.
Finally, $\langle \QQ \rangle$ denotes the ideal of $\FF[[x]]$ made of
elements with coefficients $\frac{c}{c'}$ such that $c\in \QQ$. Remark
that the latter is a priori different from $\FF[[x]] \QQ$ since we don't
suppose $\CC$ to be noetherian.

Now for an element $f$ in $\CC[[x]] \smallsetminus \QQ[[x]]$ or more
generally in $\CC_\QQ[[x]] \smallsetminus \langle \QQ \rangle$,
let us write $f=\sum_\alpha \frac{c_\alpha}{c'_\alpha} x^\alpha$ with
$c'_\alpha \in \CC \smallsetminus \QQ$. Then denote by $\exp^\modQ(f)$
the maximum (w.r.t. $\prec$) of the $x^\alpha$ such that $c_\alpha
\notin \QQ$. This is the leading exponent of $f$ modulo $\QQ$.
In the same way, we define the leading term $\lt^\modQ(f)$, leading
coefficient $\lc^\modQ(f)$ and leading monomial $\lm^\modQ(f)$ modulo $\QQ$.

If $f,g \in \FF[[x]]$ then $\exp^\modQ(fg)= \exp^\modQ(f)+\exp^\modQ(g)$
as for the usual leading exponent. However, there are some differences
with the usual situation, for example the leading coefficient $\modQ$
of $fg$ is not equal to the product of that of them. They are equal
only modulo $\QQ$ so we will have to be careful.

Now for an ideal $J\subset \CC[[x]]\smallsetminus \QQ[[x]]$, we define:
$\Exp^\modQ(J)=\{\exp^\modQ(f) | f\in J\}$.
This is a subset of $\N^n$ which is stable by sums. Thus by Dickson lemma:
\begin{equation}\label{eq:gensb}
\exists \{g_1,\ldots,g_r\}\subset J \text{ such that } \Exp^\modQ(J)=
\bigcup_j (\exp^\modQ(g_j)+ \N^n).
\end{equation}
This shall be a generic standard basis of $J$ on $V(\QQ)$.
However, this is not the definition we will adopt. In fact in the next
paragraph we will define the notion of generic reduced standard basis and
it will not be in the ring $\CC[[x]]$ so we need a more general definition:

\begin{defin}\label{def:genSB}
A generic standard basis (gen.s.b for short) of $J$ on $V(\QQ)$ is a
couple $(\GG,h)$ where
\begin{itemize}
\item[(a)]
$h\in \CC \smallsetminus \QQ$,
\item[(b)]
$\GG$ is a finite set in the ideal $\CC[h^{-1}][[x]] \cdot J$ and for any
$g\in \GG$ the numerator of $\lc^\modQ(g)$ divides $h$,
\item[(c)]
$\dps \Exp^\modQ(J)=\bigcup_{g\in \GG}(\exp^\modQ(g) +\N^n)$.
\end{itemize}
\end{defin}

Above in (\ref{eq:gensb}), $(\{g_1,\ldots,g_r\}, \prod_j \lc^\modQ(g_j))$
is a gen.s.b of $J$ on $V(\QQ)$.\\
Remark that another way to state (b) is: For any $\PP \in V(\QQ)
\smallsetminus V(h)$, the specialization $(g)_\PP$ is well defined and
belongs to $(J)_\PP$ and $\exp((g)_\PP)$ is equal to $\exp^\modQ(g)$. Remark
that $V(\QQ) \smallsetminus V(h)$ is non empty since $h\notin \QQ$.

\begin{prop}[Division modulo $\QQ$]\label{prop:divmodQ}
Let $h\in \CC \smallsetminus \QQ$ and $g_1,\ldots,g_r$ be in $\CC[h^{-1}][[x]]
\smallsetminus \QQ[h^{-1}][[x]]$ such that each $\lc(g_j)$ divides $h$.
Let $\Delta_1 \cup \cdots\cup \Delta_r \cup \bar{\Delta}$ be the partition
of $\N^n$ associated with the $\exp^\modQ(g_j)$.
Then for any $f$ in $\CC[h^{-1}][[x]]$, there exist $q_1,\ldots,q_r, R, T
\in \FF[[x]]$ such that
\begin{itemize}
\item[(o)]
$f=\sum_j q_j g_j +R+T$,
\item[(i)]
$\ND(q_j)+\exp^\modQ(g_j) \subset \Delta_j$ if $q_j\ne 0$,
\item[(ii)]
$\ND(R) \subset \bar{\Delta}$,
\item[(iii)]
the $q_j$ and $R$ are in $\CC[h^{-1}][[x]]$ and $T$ is in $\QQ[h^{-1}][[x]]$.
\end{itemize}
Moreover, $(q_1,\ldots,q_r,R)$ is unique modulo $\QQ[h^{-1}][[x]]$.
\end{prop}

\begin{proof}
Write $g_j=g_j^{(1)} - g_j^{(2)}$ with $g_j^{(2)} \in \langle \QQ \rangle$
and $\exp (g_j^{(1)})=\exp^\modQ(g_j)$ then divide $f$ by the $g_j^{(1)}$'s
in $\FF[[x]]$ as in theorem \ref{theo:div}: $f=\sum_j q_j g_j^{(1)} +R$.
We have, $f= \sum_j q_j g_j +R +T$ with $T=\sum_j q_j g_j^{(2)}$.
Conditions (i) and (ii) are satisfied by theorem \ref{theo:div}.
The third one is a direct consequence of lemma \ref{lem:denom}.
Let us prove the last statement. Since the division
takes place in $\CC[h^{-1}][[x]]$, we can specialize into $\QQ$:
$(f)_\QQ=\sum_j (q_j)_\QQ (g_j)_\QQ +(R)_\QQ$. Statements (i) and (ii) become
$\ND((q_j)_\QQ)+\exp^\modQ((g_j)_\QQ) \subset \Delta_j$ and $\ND((R)_\QQ)
\subset \bar{\Delta}$. Since $\exp^\modQ(g_j)=\exp((g_j)_\QQ)$, the latter
is exactly the result of the division of $(f)_\QQ$ by the $(g_j)_\QQ$
in $\FF(\QQ)[[x]]$ but in this division the quotients and the remainder
are unique. This implies our desired statement.
\end{proof}

\begin{cor}\label{cor:prop:divmodQ}
Retain the hypotheses and the notations of the previous proposition.
For any $\PP \in V(\QQ) \smallsetminus V(h)$, the division of $(f)_\PP$
by the $(g_j)_\PP$ is $(f)_\PP=\sum_j (q_j)_\PP (g_j)_\PP +(R)_\PP$.
\end{cor}

The proof is left to the reader (it uses the same arguments as above).
Now, here is a result similar to proposition \ref{prop:rest0}.

\begin{prop}\label{prop:rest0modQ}
Let $(\GG,h)$ satisfiying condition (a) and (b) of
definition \ref{def:genSB}. Then the following statements are equivalent.
\begin{itemize}
\item[(1)]
Condition (c) of definition \ref{def:genSB} is satisfied.
\item[(2)]
For any $f\in J$, the remainder $R$ of the division modulo $\QQ$ of $f$
by $\GG$ is $0$ (modulo $\QQ$).
\item[(3)]
The specialization $(\GG)_\QQ$ is a standard basis of $(J)_\QQ$.
\end{itemize}
\end{prop}

\begin{proof}
$(1) \Rightarrow (2)$: Let $f$ be in $J$, set $\GG=\{g_1,\ldots, g_r\}$ and
take the notations of the previous theorem, so we make the (usual) division
of $f$ by the $g_j^{(1)}$ and we want to prove that the remainder is zero.
As we recalled in subsection \ref{subsec:div}, this gives rise to sequences
 $q_j^{(i)}$, $f^{(i)}$ and $R^{(i)}$ in $\FF[[x]]$ such that for any $i$,
\[f=\sum_j q_j^{(i)} g_j^{(1)} + R^{(i)} + f^{(i)}.\]
Let us prove by an induction on $i$ that for any $i\ge 0$, the
following holds:
\[
\begin{cases}
(1,i) \quad h^l f^{(i)} \in J +\QQ[[x]] \text{ for some } l \in \N\\
(2,i) \quad R^{(i)}=0.
\end{cases}
\]
Those two statements are true for $i=0$. Assume $(1,i)$ is true then
by hypothesis (1), $\exp^\modQ(f^{(i)}) \in \exp^\modQ(g_j)+ \N^n
=\exp(g_j^{(1)})+\N^n$ for some $j$ thus $R^{(i+1)}=R^{(i)}$ which is zero by
$(2,i)$, thus $(2,i+1)$ is true. Therefore, $f^{(i+1)}=f^{(i)}-
\frac{\lc(f^{(i)})}{\lc(g_j^{(1)})} \cdot x^{\exp(f^{(i)})-\exp(g_j^{(1)})}
\cdot g_j^{(1)}$ which implies $(1,i+1)$.
The induction is done. It follows that $R=0$ and we are done.\\
$(2) \Rightarrow (1)$ is easy and left to the reader.\\
$(3) \Rightarrow (2)$: Let $f$ be in $J$ and consider the division $\modQ$
of $f$ by $\GG$: $f=\sum_j q_j g_j +R+T$. By the corollary above applied
to $\PP=\QQ$, the equality $(f)_\QQ=\sum_j (q_j)_\QQ (g_j)_\QQ +(R)_\QQ$
is the result of the division of $(f)_\QQ$ by $(\GG)_\QQ$. But $f\in J$
so $(f)_\QQ\in (J)_\QQ$ so by assumption $(R)_\QQ=0$ which implies $(2)$.\\
$(2) \Rightarrow (3)$ follows from the next theorem applied to
$\PP=\QQ$.
\end{proof}

The following lemma will be useful for our construction of comprehensive
standard bases.

\begin{lem}\label{lem:genSBinJ}
Let $\ring_0$ be a subring of $\CC[[x]]$ such that $\CC[x] \subset \ring_0$.
Let $J$ be an ideal in $\ring_0$ and $\hat{J}$ be its extension in $\CC[[x]]$.
Then there exists a gen.s.b $(\GG,\ \prod_{g\in \GG} \lc^\modQ(g))$ of
$\hat{J}$ on $V(\QQ)$ such that $\GG \subset J$.
\end{lem}
\begin{proof}
By definition and noetherianity of $(\hat{J})_\QQ$,  there exists
$f_1,\ldots,f_q \in J$ such that the $(f_j)_\QQ$ generates $(\hat{J})_\QQ$.
So by lemma \ref{lem:algebraicSB}, there exists a standard basis $G$ of
$(\hat{J})_\QQ$ which is included in $\sum_j \FF(\QQ)[x] (f_j)_\QQ$.
Therefore by multiplying each $g\in G$ by some element in $\FF(\QQ)
\smallsetminus (0)$, we may assume that each $g\in G$ is equal to
some $(f)_\QQ$ with $f\in J$.
By lifting from $\FF(\QQ)$ to $\FF$, we obtain a set $\GG \subset J$ such
that $(\GG)_\QQ=G$ is a standard basis of $(\hat{J})_\QQ$. If we define $h$
as the product of the leading coefficients $\modQ$ of $g\in \GG$, then
$(\GG,h)$ statisfies statement $(3)$ of Prop.~\ref{prop:rest0modQ}.
\end{proof}

The main result concerning generic standard basis is the following.

\begin{theo}\label{theo:genSB}
Let $(\GG,h)$ be a gen.s.b of $J \subset \CC[[x]]$ on $V(\QQ)$.
Then for any $\PP \in V(\QQ) \smallsetminus V(h)$:
\begin{itemize}
\item[(i)] $(\GG)_\PP \subset (J)_\PP$,
\item[(ii)] $\dps \Exp((J)_\PP)=\bigcup_{g\in \GG} (\exp((g)_\PP)+\N^n)
=\Exp^\modQ(J)$.
\end{itemize}
\end{theo}
In other words, $(\GG)_\PP$ is a standard basis of $(J)_\PP$ for a generic
$\PP \in V(\QQ)$ and $\Exp((J)_\PP)$ is (generically) constant and equal
to $\Exp^\modQ(J)$.

\begin{proof}
Statement (i) follows from \ref{def:genSB} (b).
Let us prove (ii). The second equality is straightforward. To prove the
first one, we shall use the criterion involving the S-functions.
Set $\GG=\{g_1,\ldots,g_r\}$ and let $f$ be in $J$ then by statement (2)
in the previous proposition, the remainder of the division modulo $\QQ$ of
$f$ by $\GG$ is $0$: $f=\sum_j q_j g_j +T$ with $T\in \langle \QQ \rangle$.
Thus by corollary \ref{cor:prop:divmodQ}, $(f)_\PP=\sum_j (q_j)_\PP
(g_j)_\PP$. As a consequence, $(\GG)_\PP$ generates $(J)_\PP$
over $\FF(\PP)[[x]]$.

Now, let us fix $g$ and $g'$ in $\GG$.
Put $S=\lc^\modQ(g') m g - \lc^\modQ(g) m'g'$ where $m=m_0/\lt^\modQ(g)$,
$m'=m_0/\lt^\modQ(g')$ and $m_0=\mathrm{lcm}(\lt^\modQ(g), \lt^\modQ(g'))$.
Remark that for any $\PP \in V(\QQ)$ such that $h \notin \PP$,
we have $(S)_\PP= S((g)_\PP,(g')_\PP)$.
Consider the division modulo $\QQ$ of $S$ by the $g_j$'s for which the
remainder $R$ is zero: $S=\sum_{j=1}^r q_j g_j +T$. By corollary
\ref{cor:prop:divmodQ}, the division of $(S)_\PP=S((g)_\PP, (g')_\PP)$
by $(\GG)_\PP$ has a zero remainder. We conclude with proposition
\ref{prop:Scriterion}.
\end{proof}

\subsection{Generic reduced standard bases}

The next result shall concern the existence of \emph{the} reduced
generic standard basis on $V(\QQ)$ (in fact we shall see that it is
unique ``modulo $\QQ$''). The importance of reduced standard bases
is well known. Reduced generic standard bases are also important.
For example, they played a fundamental role in our study of parametric
Gr\"obner fans \cite{B-genGF}.\\
Moreover, generic reduced standard bases constitute the main difference
between ``global'' and ``local'' situations. Indeed, if we study parametric
Gr\"obner bases for ideals in $\k[a][x]$ with a well-order $\prec$ then
generic reduced Gr\"obner bases have denominators with bounded
multiplicities but as we saw in Example \ref{exemple} it is not the case
when the order is local.

Let $J$ be an ideal in $\CC[[x]]$ and $\QQ$ be a prime ideal in
$\CC$.

\begin{theo}[Def{}inition-Theorem]\label{theo:gen.red.sb}
\noindent
\begin{itemize}
\item
There exists a gen.s.b $(\GG, h)$ of $J$ on $V(\QQ)$ such that $(\GG)_\QQ$
is the reduced standard basis of $(J)_\QQ$. Such a $(\GG,h)$ is called a
generic reduced standard basis (gen.red.s.b) of $J$ on $V(\QQ)$.
\item
If $(\GG,h)$ is a gen.red.s.b on $V(\QQ)$ then 
for any $\PP \in V(\QQ) \smallsetminus V(h)$, $(\GG)_\PP$ is the reduced
standard basis of $(J)_\PP$.
\end{itemize}
\end{theo}

\noindent
Such a gen.red.s.b is unique ``modulo $\langle \QQ \rangle$''.
More precisely:
\begin{lem}
Let $(\GG,h)$ and $(\GG',h')$ be two gen.red.s.b of $J$ on $V(\QQ)$
then
\begin{itemize}
\item
their cardinality and the set of their leading exponents $\modQ$ are
equal,
\item
if $g\in \GG$ and $g'\in \GG'$ satisfy $\exp^\modQ(g)=\exp^\modQ(g')$
then $g-g'$ belongs to $\QQ[(hh')^{-1}][[x]]$.
\end{itemize}
\end{lem}

\begin{proof}
The first statement is trivial by unicity of reduced standard bases.
For the second one, we have $g-g' \in \CC[hh'^{-1}][[x]]$ and by the same
argument of unicity, $(g)_\QQ -(g')_\QQ=0$ thus $g-g'\in \QQ[hh'^{-1}][[x]]$.
\end{proof}

\begin{proof}[Proof of the Theorem]
For the first statement, let $(\GG_0, h)$ be any gen.s.b of $J$ on $V(\QQ)$.
Set $\GG_0=\{g_1,\ldots, g_r\}$. By removing the unecessary elements, we may
assume that it is minimal. For any $j$ we may assume $\lc^\modQ(g_j)$
to be unitary. For any $j$, let $r_j$ be the remainder $\modQ$ of the
division modulo $\QQ$ of $g_j-\lm^\modQ(g_j)$ by $\GG_0$.
Set $\GG=\{\lm^\modQ(g_j)+r_j | j=1,\ldots,r\}$. It is easy to check that
$(\GG,h)$ is a red.gen.s.b.

Let us prove the second statement. Let $(\GG,h)$ be a gen.red.s.b.
First, we know that for any $\PP \in V(\QQ) \smallsetminus V(h)$,
$(\GG)_\PP$ is a standard basis of $(J)_\PP$.
Moreover it is minimal since $\Exp((J)_\PP)=\Exp((J)_\QQ)$ and
$\exp((g)_\PP)=\exp^\modQ(g)=\exp((g)_\QQ)$ for any $g\in \GG$.
The latter also implies that it is unitary. It just remains to prove
that it is reduced.
But this follows from the fact that $(\GG)_\QQ$ is reduced and that
for any $g\in \GG$, $\ND((g)_\PP) \subset \ND((g)_\QQ)$ (since $\QQ
\subset \PP$).
\end{proof}

\section{Polynomial case: an algorithmic construction}

Let us consider an ideal $J$ in $\k[a][x]$ where the system of
variables $a=(a^1,\ldots,a^m)$ is seen as a parameter, i.e. $\CC=\k[a]$.
In this section, we fix $\prec$ to be a monomial order on the terms
$x^\alpha$ but $\prec$ is taken arbitrary, i.e. we don't
suppose it to be local. The goal of this section is, given a prime
ideal $\QQ$ in $\k[a]$, to present a (finite) algorithm for computing
a generic standard basis of $J$ on $V(\QQ)$.

In order to make precise statements, let us give a definition of a
$\prec$-standard basis for an ideal $I \subset \k[x]$.

\begin{defin}\label{def:BSpoly}
A set $G \subset I$ is called a $\prec$-standard basis of $I$ if for
any $f\in I$, we have a standard representation: $f=\sum_{g\in G} q_g
g$ where $q_g \in \k[x]$ and either $q_g=0$ or $\exp_\prec(f) \succeq
\exp(q_g g)$.
\end{defin}

In the litterature we can find other (non equivalent) definitions
depending on the situation. When $\prec$ is a well order, this
definition coincides with:
\begin{equation}\label{eq:defSB}
G \subset I \textrm{ and } \Exp_\prec(I)=
\bigcup_{g\in G}(\exp_\prec(g) +\N^n).
\end{equation}

\begin{rem}
\begin{itemize}
\item
In general, definition \ref{def:BSpoly} only implies (\ref{eq:defSB}).
\item
If $\prec$ is local then it is well known that if $G$ satifies
(\ref{eq:defSB}) then it is a $\prec$-standard basis of $\hat{I}=\k[[x]] I$.
\end{itemize}
\end{rem}

\begin{proof}[Sketch of proof of the second statement]
By Robbiano's theorem $\prec$ is equivalent to some $<_w$ with $w\in
\R_{\le 0}^n$ (see the notations after Th. \ref{theo:div}). By a small
pertubation, we may assume $w$ to have non zero coefficients. Then we
can prove the following fact: for any $g\in \hat{I}$, there exists
$f\in I$ such that $\exp_{<_w}(g)= \exp_{<_w}(f)$ (by a truncation of
$g$). This concludes the proof.
\end{proof}

We shall separate the case when the order $\prec$ is a well order from
the case when it is not.

\subsection{The order $\prec$ is a well order}

Let us fix any well order $<_0$ on the terms $a^\gamma$,
$\gamma \in \N^m$. Let us define the order $<$:

$x^\alpha a^\gamma < x^{\alpha'} a^{\gamma'} \iff
\begin{cases}
x^\alpha \prec x^{\alpha'}\\
\text{or equality and } a^\gamma <_0 a^{\gamma'}.
\end{cases}$

Note that this order is also a well order. So in the following we will
use the usual theory of Gr\"obner bases in polynomial rings \cite{buchberger}
(see \cite{clo}).

\noindent
{\bf Note.} For an element $f \in \k[a,x]$, we can consider two
types of leading exponents (and of leading terms, coefficients, etc):
$\exp_{\prec}(f) \in \N^{n}$ and $\exp_{<}(f)\in \N^{m+n}$. Thus
we will have to be careful.

\begin{rem}\label{rem:exp_<}
For any $f\in \k[a,x]$, $\exp_<(f)=(\exp_\prec(f), \exp_{<_0}(\lc_\prec(f)))$.
\end{rem}

\

Let $G$ be a minimal Gr\"obner basis of $\tilde{J}=J+ \k[a,x]
\cdot \QQ$ w.r.t. $<$.
\begin{lem}
For any $g\in G$: $\lc_\prec(g)\in \QQ \iff g\in \QQ$.
\end{lem}
\begin{proof}
The right-left implication is trivial. Let us prove the converse one.
Write $g=\lc_\prec(g) x^\alpha+ \cdots$ where $\alpha=\exp_\prec(g)$.
We have $\exp_<(g)=\exp_<(\lc_\prec(g)) + (\alpha,0)$ but since $\lc_\prec(g)
\in \QQ \subset \tilde{J}$, we have $\exp_<(\lc_\prec(g)) \in \exp_<(g')
+(\alpha', \gamma')$ for some $g'\in G$ and $(\alpha',\gamma')\in
\N^{n+m}$. By minimality of $G$, we must have $(\alpha,0)=(\alpha',
\gamma')=(0,0)$. Therefore $g=\lc_\prec(g)$.
\end{proof}

The following proposition is probably ``well known to specialists''.

\begin{prop}
Take $\GG=G \smallsetminus \QQ$ and let $h$ be the product of
$\lc_\prec(g)$ for $g\in \GG$ then for any $\PP \in V(\QQ) \smallsetminus
V(h)$, $(\GG)_\PP$ is a $\prec$-Gr\"obner basis of $(J)_\QQ$.
\end{prop}

In this proposition, $\GG$ is not a gen.s.b strictly speaking since $\GG$
is not necessarily in $J$. So we have to (re)construct a gen.s.b from
this $\GG$: Let us denote by  $f_1,\ldots, f_q$ and $c_1,\ldots, c_e$ the
given generators of $J$ and $\QQ$ repectively.
The Gr\"obner basis $G$ is constructed by using Buchberger algorithm
starting from $\{f_1,\ldots, f_q, c_1,\ldots,c_e\}$. This
calculation is based on division of $S$-polynomials. Now, if we keep
all these divisions, it will be possible for any $f\in G$
to write explicitely $f=\sum_j u_{f,j} \cdot f_j+\sum_j v_{f,j}
\cdot c_j$. As a consequence the set $\GG'$=$\{\sum_j u_{f,j} \cdot f_j
| f \in \GG\}$ coupled with the $h$ of the theorem form the desired
generic standard basis.

\begin{proof}[Proof of the Proposition]
It is easy to see that for any $f' \in (J)_\PP$, there exists $f\in
\tilde{J}$ with $\lc_\prec(f) \notin \QQ$ such that $\exp_\prec(f')=
\exp_\prec(f)$. We have to prove that $\exp_\prec(f)\in \exp_\prec(g)+\N^n$
for some $g \in \GG$. By construction $<$ is an elimination order for
the variables $x_i$ thus $G \cap \k[a]$ is a $<_0$-Gr\"obner
basis of $\tilde{J} \cap \k[a]$. Let us treat two cases.\\
$\bullet$ $\QQ$ is included in but different from $\tilde{J} \cap \k[a]$:
in this case, there exists $g \in \k[a] \smallsetminus \QQ$ in $\GG$ and
the conclusion is trivial.\\
$\bullet$ $\QQ=\tilde{J} \cap \k[a]$: In this case, $G \cap \QQ=G \cap \k[a]$
and this set is $<_0$-Gr\"obner basis of $\QQ$. Write $f=\lc_\prec(f)
\lt_\prec(f)+ f_0$ and divide $\lc_\prec(f)$ by $G \cap \QQ$ w.r.t. $<_0$.
Let $r$ be the remainder. It is not zero. Put $f_1= r\cdot \lt_\prec(f)+
f_0$. We have $\exp_\prec(f)=\exp_\prec(f_1)$ and $\lc_\prec(f_1)=r \notin
\QQ$. By construction $f_1 \in \tilde{J}$ so there exists $g\in G$ such that
$\exp_<(f_1) \in \exp_<(g)+\N^{n+m}$. By Remark \ref{rem:exp_<}, this
implies $\exp_\prec(f_1)\in \exp_\prec(g)+\N^n$ and $\exp_{<_0}(r) \in
\exp_{<_0}(\lc_\prec(g))+\N^m$.
Suppose $g\in \QQ$ then by the lemma above, $g=\lc_\prec(g)\in \QQ$
and $\exp_{<_0}(r) \in \exp_{<_0}(g) +\N^m$ but this is impossible
since $r\notin \QQ$. Thus $g$ must belong to $\GG$ and we are done.
\end{proof}

\subsection{The order $\prec$ is not a well order}

In this case, our method is based on a homogenization and a computation
w.r.t. to a well order (as introduced by D.~Lazard \cite{lazard}).

As above, let us fix any well order $<_0$ on the terms $a^\gamma$,
$\gamma \in \N^m$. Let $z$ be a new variable and let us define the
orders $\prec^z$ and $<^z$:

$x^\alpha z^k \prec^z x^{\alpha'} z^{k'} \iff
\begin{cases}
|\alpha|+k < |\alpha'|+k'\\
\text{or equality and } x^\alpha \prec x^{\alpha'},
\end{cases}$

$a^\gamma x^\alpha z^k <^z a^{\gamma'} x^{\alpha'} z^{k'} \iff
\begin{cases}
x^\alpha z^k \prec^z x^{\alpha'} z^{k'}\\
\text{or equality and } a^\gamma <_0 a^{\gamma'}.
\end{cases}$

Here above, $\alpha$, $\gamma$, $k$ are in $\N^n$, $\N^m$ and
$\N$ respectively. Note that these orderings are well orders.

For any ring $A$ and any $f$ in $A[x]$, write $f=\sum_\alpha c_\alpha
x^\alpha$ and define the homogenization of $f$ in $A[x,z]$ as $h(f)=
\sum_\alpha c_\alpha x^\alpha z^{d-|\alpha|}$ where $d$ is the total degree
of $f$ in the variables $x_i$. More generally an element $f$ of the form
$f=\sum_{\alpha,k} c_{\alpha,k} x^\alpha z^k$ with $c_{\alpha,k}=0$
if $|\alpha|+k\ne d$ is called homogeneous (of degree $d$).\\

Now let $f_j$, $j=1,\ldots, q$ be the given generators of $J \subset \k[a,x]$
and set $h(J) \subset \k[x,a,z]$ to be the ideal generated by the $h(f_j)$.

\begin{theo}\label{theo:genSBpoly}
Let $G$ be a $<^z$-standard basis of $J'=h(J)+ \k[x,a,z]\cdot
\QQ$ made of homogeneous elements. Then the set $\GG=\{f_{|h=1};
f\in G \smallsetminus \QQ\}$ satisfies the following: Let $h$ be the product
of the $\lc_\prec(g)$, $g\in \GG$, then there exists $h' \in \k[a]
\smallsetminus \QQ$ such that for any $\PP \in V(\QQ) \smallsetminus V(hh')$,
$(\GG)_\PP$ is a $\prec$-standard basis of $(J)_\PP$ following definition
\ref{def:BSpoly}.
\end{theo}
\noindent
{\bf Remarks.} (1) Since $J'$ is generated by homogeneous elements and
Buchberger algorithm conserves homogeneity, it is always possible to
construct $G$ made of homogeneous elements.
(2) As in the previous subsection, $\GG$ is not strictly speaking a gen.s.b
and we can reconstruct a gen.s.b from $G$.

\begin{lem}
Take $\PP \in \spec(\k[a])$ then
\begin{itemize}
\item[(a)]
$(J)_\PP=\sum_{j=1}^q \FF(\PP)[x] \cdot (f_j)_\PP$.
\item[(b)]
There exists $h' \in \k[a] \smallsetminus \QQ$ such that if $h' \notin \PP$
then $(h(J))_\PP=\sum_{j=1}^q \FF(\PP)[x,z] \cdot h((f_j)_\PP)$.
\end{itemize}
\end{lem}

\begin{proof}
Statement (a) is trivial, let us prove (b). For each $j$, let $c_j \in \k[a]$
be a coefficient of some term of $f_j \in \k[a][x]$ with maximal degree.
Put $h'=\prod_j c_j$. For any $\PP \notin V(h')$, the degree of $f_j$ is
equal to that of $(f_j)_\PP$ so $(h(f_j))_\PP=h((f_j)_\PP)$ from
which the statement follows.
\end{proof}

\begin{lem}
For any $\PP \in V(\QQ) \smallsetminus V(hh')$, $(G \smallsetminus \QQ)_\PP$
is a $\prec^z$-standard basis of $\sum_{j=1}^q \FF(\PP)[x,z] \cdot
h((f_j)_\PP)$ and is made of homogeneous elements.
\end{lem}

\begin{proof}
By the previous lemma, since $\PP \in V(\QQ) \smallsetminus
V(h')$, $(G)_\PP$ generates the ideal in question.
We conclude with the proposition above.
\end{proof}

The following lemma is a classical result.
\begin{lem}\label{lem:classic}
Let $I$ be an ideal in $\k[x]$ generated by $f_1,\ldots,f_q$.
Let $G$ be a homogeneous $\prec^z$-standard basis of the ideal $h(I)$
generated by $h(f_j)$, $j=1,\ldots,q$. Then $G_{|z=1}$ is a
$\prec$-standard basis of $I$ in the sense of definition
\ref{def:BSpoly}.
\end{lem}
\begin{proof}
Let $f\in I$. Write $f=\sum_j u_j f_j$. Homogenization implies that there
exist $l, l_1,\ldots,l_q \in \N$ such that $z^l h(f)=\sum_j z^{l_j}
h(u_j)h(f_j)$ so $z^l h(f)$ belongs to $h(I)$. By definition of $G$:
$z^l h(f)=\sum_j q_j g_j$ where $G=\{g_1,\ldots, g_r\}$ and $q_j \in
\k[x,z]$ and $\exp_{\prec^z}(f)\ \succeq^z \exp_{\prec^z}(q_j g_j)$.
By division, the $q_j$ are homogeneous. But for a homogeneous element
$H\in \k[x,z]$, we have $\pi(\exp_{\prec^z}(H))=\exp_\prec(H_{|z=1})$
where $\pi(\alpha,k)=\alpha$. Thus specializing $z=1$ gives the desired
standard representation.
\end{proof}

Now the proof of Theorem \ref{theo:genSBpoly} is a direct application
of this lemma to our situation.

\section{Illustration}

As we said, in another paper \cite{B-polygen} we use generic standard bases
for studying the local Bernstein polynomial for a deformation of a
singularity. However, in order to keep a reasonable size to the present paper
we shall restrict ourselves to two direct applications.

\subsection{Comprehensive standard bases}

Our goal here is not to give a general theory of comprehensive standard
bases (which in the global case, were treated by V.~Weispfenning in
\cite{weisp92}, \cite{weisp03}), we only intend to illustrate in a natural
situation how we can use generic standard bases.

Let $A \subset \C^m$ and $X \subset \C^n$ be polydisks centered at $0$ and
let $J$ be an ideal in $\mathcal{O}_{A\times X}$ which denotes the ring
of analytic functions on $A\times X$.
For $a_0 \in A$, we denote by $J_{|a_0} \subset \mathcal{O}_X$ the ideal
obtained by specializing $a=a_0$. This ideal can be identified with the
specialization $(J)_{m_{a_0}}$ where $m_{a_0} \subset \mathcal{O}_A$ is
the maximal ideal generated by the $a^i-a_0^i$, $i=1,\ldots,m$.

For $Y \subset A$, a subset $W \subset Y$ is locally closed if it is the
difference of two (analytic) closed subsets of $Y$. $W$ is constructible
if it is a finite union of locally closed subsets.
\begin{prop}\label{prop:comprehensive}
\noindent
\begin{itemize}
\item[(i)]
There exists a finite partition $A=\cup_k W_k$ into analytic locally closed
subsets of $A$ such that for any $k$, $\Exp(J_{|a_0})$ is constant
for $a_0 \in W_k$.
\item[(ii)]
There exists finite set $\GG \subset J$ such that for any $a_0 \in A$,
$\GG_{|a_0}$ is a standard basis of $J_{|a_0}$.
\end{itemize}
\end{prop}

\begin{proof}
We shall prove both statements at the same time.
By induction on the dimension, let us prove that for any Zariski closed
set $Y \subset A$, (i) and (ii) where we replace $A$ by $Y$ are true.
If $\dim Y=0$, $Y$ is a finite union of points so (i) is trivial and for
(ii), one has to take a generic standard basis of $\O_A[[x]] \cdot J$ w.r.t
the maximal ideals associated to each point. Thanks to lemma
\ref{lem:genSBinJ}, it can be chosen in $J$.

Suppose $\dim Y \ge 1$. Write $Y=V(\QQ_1) \cup \cdots \cup V(\QQ_r)$
(here $V(\QQ_i)$ is the ``usual'' zero set in $A$) with $\QQ_i$ prime
in $\mathcal{O}_A$. For each $Q_i$, let $(\GG_i, h_i)$ be a generic
standard basis of the extension $\O_A[[x]]\cdot J$. It can be chosen in $J$.
Now write $Y=Y_1 \cup Y_2$ where
$Y_1=\bigcup_i (V(\QQ_i) \smallsetminus V(h_i))$ and $Y_2=\bigcup_i
(V(\QQ_i) \cap V(h_i))$. Set $\GG'=\cup_i \GG_i$. For any $a_0$ in $Y_1$,
$\GG'_{|a_0}$ is a standard basis of $J_{|a_0}$.
We have $\dim Y_2 < \dim Y$ so let us apply the induction hypothesis to $Y_2$:
we obtain a finite set $\GG'' \subset J$ such that for any $a_0 \in Y_2$,
$\GG''_{|a_0} \subset J_{|a_0}$ is a standard basis; we also obtain that
$Y_2$ is a finite union of locally closed sets such that on each of them
the map $a_0 \mapsto \Exp(J_{|a_0})$ is constant.
Finally we set $\GG=\GG' \cup \GG''$ and we reorganize the writing of $Y$
in order to have a partition (recall that constructible sets are stable
by intersection, finite union and complementation).
\end{proof}

\begin{rem}
Suppose the order $\prec$ is not necessarily local and let $I$ be an ideal
in $\k[a][x]$. Take the notations of the previous section.
If $\GG$ is a homogeneous comprehensive Gr\"obner basis of $h(I)$ for
$\prec^z$ (with the definition of \cite{weisp92}) then by lemma
\ref{lem:classic}, $\GG_{|z=1}$ is a comprehensive standard basis of
$I$ for $\prec$.
\end{rem}

\subsection{Hilbert polynomial}

Let $J$ be an ideal in $\mathcal{O}_{A\times X}$ (we keep the notations
above).
\begin{prop}
The partition of $A$ given by the local Hilbert polynomial of $\mathcal{O}_X
/J_{|a}$ at $x=0$ is constructible.
\end{prop}
Given an analytic function $f\in \mathcal{O}_{A\times X}$ such that
$f(0, a)=0$ for any $a\in A$. Then applying this proposition to the ideal
generated by the partial derivatives $\frac{\partial f}{\partial x_i}$
will provide a constructible partition of $A$ such that the Milnor number
of $f_{|a_0}$ is constant on each strata (this result can also be derived
from the semi-continuity of the Milnor number: see \cite{broughton} when
$f$ is polynomial).
For the definition of the local Hilbert polynomial, one can refer to
\cite{matsumura} and \cite{singular}. By an abuse of notations, we will
identify $J$ with its germ in $\C\{a,x\}$ and $J_{|a}$ with its germ in
$\C\{x\}$.

Denote by $m$ the maximal ideal in $\C\{x\}$.
For an ideal $I$ in $\C\{x\}$, the (local) Hilbert-Samuel function of $I$:
$HSF_I:\N \to \N$ is defined by $HSF_I(r)=\dim_\C (\C\{x\}/(I+m^{r+1}))$.
For a set $E \in \N^n$ such that $E+\N^n=E$, we define its Hilbert-Samuel
function $HSF_E:\N \to \N$ as $HSF_E(r)= \mathrm{card}\{\alpha\in \N^n;
\, \alpha\in \N^n \smallsetminus E, \, |\alpha|\le r\}$.

There exists a rational polynomial $HSP_I$ (the local Hilbert polynomial)
such that for $r \in \N$ large enough  $HSP_I(r)=HSF_I(r)$.

\begin{lem}
Let $\prec$ be a local order such that: $|\alpha| < |\alpha'| \Rightarrow
x^\alpha \succ x^{\alpha'}$. Set $E=\Exp_\prec(I)$ then $HSF_I=HSF_E$.
\end{lem}
The proof of this lemma is easy and left to the reader. Now the proposition
follows easily from this lemma and Prop.~\ref{prop:comprehensive}.

\section{Extension to differential operators rings}

Given an ideal $J$ in $\CC[[x]]$ or more generally in a subring $\ring$ of
$\CC[[x]]$ we have shown the existence of gen.s.b of the extension
$\hat{J}=\CC[[x]]J$ of $J$. We have seen that a gen.s.b can be chosen
in $\hat{J}$ (see (\ref{eq:gensb})), while red.gen.s.b are not in $\CC[[x]]$
in general. We have also seen that the use of truncated divisions in
Buchberger algorithm allows us to construct a gen.s.b in $J$ itself.

The only results that we needed for this construction were: a formal
division procedure, the fact that Buchberger algorithm works and
truncated divisions. This means that our construction works in many other
situations.
Let us state this construction for two of them namely: rings of
differential operators with parameters and the $(0,1)$-homogenization of
the latter. In \cite{B-polygen}, we used gen.s.b in rings of differential
operators.\\

Her we will follow (Castro-Jim\'enez, Granger, \cite{cg}).
Let $\FDn(\k)$ be the ring of differential operator with coefficients
in $\k[[x]]$ and $x=(x_1,\ldots,x_n)$. We denote by $\dx{1}, \ldots, \dx{n}$
the derivations. An element $P$ in this ring has a unique writing:
$P=\sum_{\alpha, \beta} c_{\alpha, \beta} x^\alpha \ddx^\beta$.
Define its Newton diagram $\ND(P)\subset \N^{2n}$ as the set of
$(\alpha,\beta)$ with $c_{\alpha,\beta}\ne 0$.
Following (\cite{cg}, chapter 2), let $w=(w^1,w^2)\in \R^{n+n}$ where
$w^2$ has strictly positive coefficients, and $w^1$ has non negative ones.
Define an order $<_w$ on $\N^{2n}$ (or equivalently on the monomials $x^\alpha
\xi^\beta$) first by $w^2$, then the inverse of $w^1$ and refine them
by the inverse of a well order $<_0$ (on $\N^{2n}$).
We can define the leading exponent, leading coefficient, leading term and
leading monomial w.r.t. $<_w$ of an element $P\in \FDn(\k)$ as in
the previous sections. For $P\in \FDn(\k)$, we denote by $\ord^{w^2}(P)$ the
maximum of $w^2\cdot \beta$ for $(\alpha,\beta) \in \ND(P)$. Remark that
this order gives rise to the $w^2$-Bernstein filtration $F^{w^2}_k=
\{P|\ord^{w^2}(P)\le k\}$ for which the graded ring $\gr^{w^2}(\FDn(\k))$
is isomorphic to $\k[[x]][\xi]$. We see then that the order $<_w$ is
``adapted'' to the weight vector $w^2$. Such orders are used for the
calculation of the characteristic variety of an analytic or formal $D$-module.

With the order $<_w$ we have a division theorem similar to Th.
\ref{theo:div}, where we just replace $\k[[x]]$ by $\FDn(\k)$ and $\C\{x\}$
by $\mathcal{D}_n$, and $\prec$ by $<_w$, see (\cite{cg}, Th 2.4.1).
As we did in the previous sections, truncated divisions work here again.
The same definition of standard basis gives rise to the same criterion
as in Prop. \ref{prop:rest0}: a system of generators is a standard basis
if the division of all the S-operators by this system has zero as the
remainder, see (\cite{cg}, Prop. 2.5.1). As a consequence, Buchberger
algorithm works in this situation (\cite{cg}, 2.5).

For an operator $P\in \FDn(\CC)$, we can define $\exp_{<_w}^\modQ(P)$, etc
as previously. Let $J$ be an ideal in $\FDn(\CC)$.
\begin{defin}
A generic standard basis of $J$ on $V(\QQ)$ w.r.t. $<_w$ is a couple
$(\GG,h)$ where\\
(a) $h\in \CC \smallsetminus \QQ$,\\
(b) $\GG$ is a finite set in the ideal $\FDn(\CC[h^{-1}]) \cdot J$ and
for any $g\in \GG$ the numerator of $\lc_{<_w}^\modQ(g)$ divides $h$,\\
(c) $\dps \Exp_{<_w}^\modQ(J)=\bigcup_{g\in \GG}(\exp_{<_w}^\modQ(g)
+\N^{2n})$.
\end{defin}
Then we can state a division modulo $\QQ$ as in Prop. \ref{prop:divmodQ} and
we obtain a criterion similar to Prop. \ref{prop:rest0modQ}. Finally here
are the analogues theorems to \ref{theo:genSB} and \ref{theo:gen.red.sb}.

\begin{theo}
Let $(\GG,h)$ be a gen.s.b of $J \subset \FDn(\CC)$ on $V(\QQ)$.
Then for any $\PP \in V(\QQ) \smallsetminus V(h)$:\\
$\bullet$ $(\GG)_\PP \subset (J)_\PP$,\\
$\bullet$ $\dps \Exp_{<_w}((J)_\PP)=\bigcup_{g\in \GG} (\exp_{<_w}((g)_\PP)
+\N^{2n}) = \Exp_{<_w}^\modQ(J)$.
\end{theo}

\begin{theo}[Def{}inition-Theorem]\

\noindent
$\bullet$
There exists a gen.s.b $(\GG, h)$ of $J$ on $V(\QQ)$ such that $(\GG)_\QQ$
is the reduced standard basis of $(J)_\QQ$. We call $(\GG,h)$ a gen.red.s.b
of $J$ on $V(\QQ)$.\\
$\bullet$
If $(\GG,h)$ is a gen.red.s.b on $V(\QQ)$ then 
for any $\PP \in V(\QQ) \smallsetminus V(h)$, $(\GG)_\PP$ is the reduced
standard basis of $(J)_\PP$.
\end{theo}

In the case treated above, we have worked with an order adapted to the
weight vector $w^2$ but for some situations we need standard bases w.r.t.
``any'' weight vector $w$. Thus as in Lazard \cite{lazard} for the
polynomial case, we have a homogenized ring which unables us to work
with any admissible weight vector $w$ (see the definition below).\\

In the following, we follow Assi et al \cite{acg}.
Let $t$ be a new variable. Define $\FDn(\k)\langle t \rangle$
as the $\k$-algebra generated by $\k[[x]]$, the $\dx{i}$'s and $t$
where the only non trivial commutation relations are $\dx{i} a(x)-
a(x) \dx{i}= \frac{\partial a}{\partial x_i} \cdot t$. If $\k=\C$ and
if we replace $\k[[x]]$ by $\C\{x\}$, we obtain $\mathcal{D}_n \langle t
\rangle$. A weight vector $w\in \R^{2n}$ is called admissible if
$w_i\le 0$ and $w_i+w_{n+i}\ge 0$ for $i=1,\ldots, n$. Given an
admissible weight vector, we can define an order $<^h_w$ on $\N^{2n+1}$
or equivalently on the monomials $x^\alpha \xi^\beta t^k$ as follows.
We define $<_w^h$ in a lexicographical way by $[|\beta|+k, \ w, \ |\beta|,
\ >_0]$ where $<_0$ is a fixed total well order on $\N^{2n+1}$.
With this order, the authors of \cite{acg} proved a division theorem
in $\mathcal{D}_n \langle t \rangle$ and in $\FDn(\k)\langle t \rangle$
for homogeneous operators as in Th.~\ref{theo:div}.
This unables us to construct generic standard bases for homogeneous ideals
in $\FD(\CC)\langle t \rangle$ w.r.t. $<_w^h$. We thus obtain the analogues
results to the two previous theorems.

\end{document}